\documentclass[twoside,11pt,reqno]{amsart}
\usepackage{amsfonts,amsmath,amssymb,amscd,mathrsfs}

\usepackage{latexsym,amssymb}
\usepackage{enumerate, verbatim, mathrsfs}
\usepackage{graphics}
\usepackage{amsthm}
\usepackage{latexsym}
\usepackage{longtable}
\usepackage{epsfig}
\usepackage{amsmath}
\usepackage{mathdots}
\usepackage{hhline}
\newtheorem{thm}{\sc Theorem}[section]
\newtheorem{prop}[thm]{\sc Proposition}
\newtheorem{lem}[thm]{\sc Lemma}
\newtheorem{cor}[thm]{\sc Corollary}

\newcommand{\pf}{{\it Proof:\quad}}
\newcommand{\dop}{\dot\partial}
\newcommand{\de}{\dot\epsilon}
\newcommand{\p}{\partial}
\newcommand{\e}{\epsilon}
\newcommand{\gl}{\lambda}
\newcommand{\gnp}{\nu^{+}}
\newcommand{\gnm}{\nu^{-}}
\newcommand{\dgnp}{\dot\nu^{+}}
\newcommand{\dgnm}{\dot\nu^{-}}

\newcommand{\gab}{\bar{\alpha}}
\newcommand{\ga}{\alpha}
\newcommand{\gbb}{\bar{\beta}}
\newcommand{\gb}{\beta}
\newcommand{\ggab}{\bar{\gamma}}
\newcommand{\gga}{\gamma}
\newcommand{\s}{\sigma}
\newcommand{\Hy}{\mathcal{H}}  
\newcommand{\F}{\mathbb{F}}

\newcommand{\dne}{\hfill $\Box$ \vspace{0.3cm}}

 \newcommand{\ssapp}{\p\!\!\downarrow\!\uparrow\!\!\e}

\begin{document}

\title[The Permutation modules of $\Hy^{n}$]{Permutation Modules associated to the \\ Hyperoctahedron
 and Group Actions}
\author{J. Siemons}
\address{School of Mathematics, University of East Anglia, Norwich, UK}
\email{j.siemons@uea.ac.uk}
\author{ B. Summers}
\address{School of Mathematics, University of East Anglia, Norwich, UK}
\email{jbensums@gmail.com}

\date{\scriptsize Version August 2018, printed \today \\[4pt]{\sc Keywords:} Hyperoctahedral complex, Permutation module, Decomposition, Orbit number and Orbit module\\[4pt]
}
\maketitle

\noindent  {\small {\sc Abstract:}\, We investigate the permutation modules associated to the set of $k$-dimensional faces of the  hyperoctahedron in dimension $n$, denoted  $\Hy^{n}.$ For any $k\leq n$ such a module can be defined over any field $\F$, it is called a {\it face module} of $\Hy^{n}$ over $\F.$ 
We describe a spectral decomposition of such face modules into submodules and show that these submodules are irreducible under the hyperoctahedral group $B_{n}.$  The same method can be used to  describe the exact relationship between the face modules in any two dimensions  $0\leq t\leq k\leq n.$   Applications of this technique include a rank formula for the rank of the incidence matrix of $t$-dimensional  versus $k$-dimensional faces of $\Hy^{n}$ and a characterization of $(t,k,\ell)$-designs on $\Hy^{n}.$ We also prove an orbit  theorem for subgroups of the  hyperoctahedral group on the set of faces of $\Hy^{n}.$ The decomposition method is elementary, mostly characteristic free and does not use the representation theory of  automorphism groups. It is therefore quite general and can be used to decompose permutation modules  associated to other  geometries. }

\section{Introduction}

In this paper we investigate  the 
permutation modules associated to the the set of faces of the $n$-dimensional cross-polytope $\Hy^{n}$, also  known as the dual of the $n$-cube or the hyperoctahedron. We regard $\Hy^{n}$ as a simplicial complex on a set $V$ of size $2n,$ the vertices of $\Hy^{n},$  and identify its $k$-dimensional faces with certain $k$-element subsets of $V.$ (See also Stanley's book~\cite{Stanley}. Our terminology departs a little from the usual one: a $k$-dimensional face of the simplicial complex has geometric dimension $k-1.$ All definitions are given in Section 2.)  While  the hyperoctahedron and its associated geometries  are among  the most widely studied objects in combinatorics  and group theory it appears that several elementary features of this simplicial complex have remained unexplored.

For the  integers $t$ and $k$ with $0\leq t\leq k\leq n$ the $(t,k)$-{\it incidence matrix} ${\mathcal S}_{t,k}$ of ${\mathcal H}^{n}$ is the $0$--$1$ matrix   whose rows and columns are indexed by the $t$- and $k$-dimensional faces of $\Hy^{n},$  respectively, so that the $(x,y)$-entry of ${\mathcal S}_{t,k}$ is $1$ if $x\subseteq y$ and $0$ otherwise.  One of the results in this paper is a formula for the rank of this incidence matrix.

\bigskip
\begin{thm} Let ${\mathcal H}^{n}$ be the $n$-dimensional hyperoctahedron  and let $0\leq t\le k\leq n.$ Then the $(t,k)$-incidence matrix of ${\mathcal H}^{n}$  has rank $${\rm rank}({\mathcal S}_{t,k})=\sum_{0\leq j\leq t}\,\min\Big\{  {n\choose t}{t\choose j},\,\,{n\choose k}{k\choose j}\Big\}$$ in characteristic $p=0$ or $p>n.$ 
\end{thm}

This result is an applications  of a new method to decompose  the permutation modules associated to  $\Hy^{n}.$ The same technique can be applied essentially to any finite partially ordered sets that has a rank function, and so is of general interest.  

\bigskip
To explain this method let $0\leq k\leq n$ and let $L_{k}^{n}$ be the set of all $k$-dimensional faces of $\Hy^{n}.$ For the field $\F$ denote the vector space over $\F$ with basis $L^{n}_{k}$ by $M^{n}_{k}=\F L^{n}_{k}.$ It  is well-known that the automorphism group of $\Hy^{n}$ is the octahedral group $B_{n}$ and so $M^{n}_{k}$ is a permutation module for $B_{n}.$ 
In Section~2\, we introduce two linear maps $$\gnp_{k}\!:\,M^{n}_{k}\to M^{n}_{k}\quad\text{and}\quad  \s_{k}\!:\,M^{n}_{k}\to M^{n}_{k}$$ which are defined  via the {containment} and the {`opposite' relation} for the faces of $\Hy^{n}.$  A key result of this paper is the following  theorem. 

\bigskip 
\begin{thm}[Spectral Decomposition]\label{Spec} Let $\F$ be a field of characteristic $p$ and  let  $M^{n}_{k}=\F L^{n}_{k}$ for $2\leq n$ and $0\leq k\leq n.$   Suppose that $p>n$ or $p=0.$ Then there are linear 
maps $\gnp_{k}$ and $\s_{k}\!:\,M^{n}_{k}\to M^{n}_{k}$ which give rise to decompositions 
\begin{eqnarray}M^{n}_{k}&=&\oplus_{j} \,E_{k,j}\,\,\,\,\text{with $0\leq j\leq k,$ and}\nonumber\\E_{k,j}&=&\oplus_{i} \,E_{k,j,i}\,\,\,\text{ with  $0\leq i\leq \min\{k-j, n-k\}$}\nonumber\end{eqnarray} 
into the eigenspaces $E_{k,j}$ of $\s_{k}$ and the eigenspaces $E_{k,j,i}$ of $\gnp_{k}$ respectively. For fixed $k$ the  $E_{k,j,i}$  are pairwise non-isomorphic irreducible $B_{n}$-modules for all $0\leq j\leq k$ and $0\leq i\leq \min\{k-j, n-k\}.$ Furthermore, $E_{k,j,i}$ is isomorphic to $E_{k',j',i'}$ if and only if  $j=j'$ and $i=i'.$
\end{thm}

We give explicit formulae for all eigenvalues  of $\gnp_{k}$ and $\s_{k},$  and  the corresponding projection maps  onto eigenspaces in Section~3 and~4\,. Theorem~\ref{Spec} follows from Theorems~\ref{fulldecomp} and~\ref{thm3.1}.

\bigskip


A particularly important role is played by Lemmas~\ref{atype}\, and ~\ref{btype} in Section 2.  These apply, with suitable  modifications, to the permutation modules associated to geometries of Type A and Type B more generally. The two lemmas  are  at the heart of the results here. We also indicate how such spectral techniques can be applied to other geometrical posets. Results of this kind will be available in forthcoming work~\cite{SZ} for symplectic and orthogonal groups.


\bigskip
In Sections~3\,  we consider other applications and mention the connection to the standard association schemes on the hyperoctahedron $\Hy^{n}$. We define $(t,k,\ell)$-designs on $\Hy^{n}$ and give a necessary and sufficient condition for a family of $k$-dimensional faces to be such a design.   In Section~4\, we prove the irreducibility of the $E_{k,j,i}$ in the main theorem and derive orbit theorems for subgroups of the hyperoctahedral group $B_{n}.$



\bigskip

\section{\sc The Hyperoctahedron  and its  associated \\ Permutation Modules}

We describe the faces of the hyperoctahedral complex and the permutation modules that are associated to this complex.  Let $n$ be a positive integer and let $V=\{\ga_{1}, \gab_{1},...,\ga_{n}, \gab_{n}\}$ be a set of $2n$ distinct elements which we call  vertices.  It is useful to put  $\bar{\gab}_{i}:=\ga_{i}.$
For $0\leq k\leq n$ let $L^{n}_{k}$ be the set of all $k$-element subsets $x$ of $V$ so that $|x\cap\{\ga_{i},\gab_{i}\}|\leq 1$ for all $1\leq i\leq n.$  Such sets $x$ are called  {\it faces.} 

We regard $L^{n}:=\bigcup_{k=0}^{n}\,L^{n}_{k}$ as a ranked partially ordered set in which the partial order is  given  by the inclusion relation for the subsets of $V.$ Geometrically $L^{n}$ is the complex of the $n$-dimensional cross polytope or hyperoctahedron on the vertex set $V.$ In order to simplify notation, a $k$-set belonging to $L^{n}_{k}$ is called a $k$-dimensional face, or just $k$-face, rather than $(k-1)$-dimensional face, as is more common for simplicial complexes. 

In this interpretation $\ga_{i}\leftrightarrow \gab_{i}$ is the pairing into {\it opposite vertices.} It naturally extends to the pairing $x=\{\gb_{1},..,\gb_{k}\}\leftrightarrow \bar{x}=\{\gbb_{1},..,\gbb_{k}\}$ for all faces of the hyperoctahedron,  associating to every face its unique {\it opposite face.}  (Here $x\subseteq \{\ga_{1}, \gab_{1},..,\ga_{n}, \gab_{n}\}$ is an arbitrary face,  having in mind that $\bar{\gab}_{i}:=\ga_{i}.)$ We refer to this structure with the given  pairing as the  {\it face complex} of the hyperoctahedron  in dimension $n$ on the vertex set $V,$ it is denoted by $\Hy^{n}.$ For any $0\leq t\leq k\leq n$ we  regard $(L^{n}_{t},\,L^{n}_{k};\,\subseteq)$ as an incidence structure in which $x\in L^{n}_{t}$ and $y\in L^{n}_{k}$ are {\it incident} with each other if and only if $x\subseteq y.$ The corresponding incidence matrix is denoted ${\mathcal S}_{t,k}.$

Let $\F$ be a field and let $X$ be a set. Then $\F X$ denotes the vector space over $\F$ with basis $X.$ The elements of $\F X$ are the formal sums $f=\sum_{x\in X}\,f_{x}\,x$ with $f_{x}\in \F.$ Naturally we identify $1x$ with $x$ so that $X\subset \F X.$  On this vector space we have the standard inner product $(-,-)$ defined by $(x,x')=1$ if $x=x'$ and $(x,x')=0$ otherwise, for all $x,\,x'\in X.$ Since $X$ is an orthonormal basis of $F X$ we have $$f=\sum_{x\in X}\,(f,x)\,x\,\,\,\text{for all}\,\,\, f\in \F X.$$ We refer to $\F X$ also as a module since it is a permutation module for any permutation group on $X,$ and when $X$ is a set of faces then $\F X$ is the face module for $X.$   

\bigskip
For the hyperoctahedron  we consider the face modules $M^{n}_{k}:=\F L^{n}_{k}$ and $M^{n}:=\F L^{n}=\bigoplus_{k=0}^{n}\, M^{n}_{k}$ in this fashion.  We have $M^{n}_{0}=\F \{\emptyset\} =\F;$ it is useful to put  $M^{n}_{k}:=0$ if $k<0$ or $k>n.$ The dimension of  $M^{n}_{k}$ is $|L^{n}_{k}|=\sum_{i=0}^{k}\,{n\choose k}{k\choose i}=2^{k}{n\choose k}.$ 
It is very useful to regard elements of $M^{n}$ as polynomials in the polynomial ring $\F[\ga_{1}, \gab_{1},..,\ga_{n}, \gab_{n}]$ and use the operations of this ring naturally. For instance, we write  $$3\{\ga,\gbb,\gga\}+2\{\ga,\gb,\epsilon\}=\ga(3\gbb\gga+2\gb\epsilon).$$ In the same spirit we suppress set-brackets and write $\ga=\{\ga\}$ etc.  where possible. (As an aside, we do not claim that $M^{n}$ is closed under this multiplication: if $x, y$ are faces then $x\cdot y$ may not be a face. However, one may put $x\cdot y=0$ in this case. This is a standard construction in the Stanley-Reisner ring of a simplicial complex, see Chapter~II in Stanley's book~\cite{Stanley}.)

\bigskip
Let $0\leq k\leq n.$ We construct a basis for $M^{n}_{k}.$ Denote $\dot{V}:=\{\ga_{1},..,\,\ga_{n}\}$ and  $\dot{L}^{n}_{k}:=\{ x\in L^{n}_{k}\,:\, x\subseteq \dot{V}\}.$ 
If  $a,\,b$ are disjoint subsets of $\dot{V}$ we put  $$[a,b]:=\Pi_{\ga\in a}(\ga+\gab)\cdot\Pi_{\gb\in b}(\gb-\gbb)\,\,\in M^{n}_{k}$$  where $|a|+|b|=k.$ In particular, $[\emptyset,\emptyset]=1\!\cdot\!\emptyset.$ For instance, when $k=2,$ $a=\{\ga\}$ and $b=\{\gb\}$ then $[a,b]=(\ga+\gab)(\gb-\gbb)=\{\ga,\gb\}  -\{\gab,\gb\} - \{\ga,\gbb\} + \{\gab,\gbb\}.$ It is useful to put $[a,b]=0$ if $a\cap b\neq \emptyset$ and to extend this notation distributively to more general expressions: If $f=f_{1}a_{1}+..+f_{t}a_{t}$ with subsets $a_{i}\subseteq \dot{V}$ and $f_{i}\in \F$ then we write $f_{1}[a_{1},b]+..+f_{t}[a_{t},b]=:[f,b].$

\bigskip
\begin{prop}\label{thm2.1} Suppose that $\F$ is a field of  characteristic  $p\neq 2$  and let $0\leq k\leq n.$ Then the elements $[a,b]$ with $a,\,b\subseteq \dot{V},$ $a\cap b=\emptyset$  and $k=|a|+|b|$ form a basis of $M^{n}_{k}.$
\end{prop}

The condition that $\F$ has characteristic $p\neq 2$ is essential.  For $p=2$ the  $[a,b]$ are no longer linearly independent. For example, we have $[\ga,\gb]=[\gb,\ga].$

\medskip 
\pf The number of pairs $(a,b)$ where $a,\,b$ are disjoint subsets of $\dot{V}$ and $|a|+|b|=k$ is $\sum_{i=0}^{k}\,{n\choose k}{k\choose i}=\dim M^{n}_{k}.$ It suffices to show that the $[a,b]$ span $M^{n}_{k}.$ This is true when $n=0$ and so we apply induction over $n.$  
Let  $n>0$ and  $x\in L^{n}_{k}.$ By induction, for $\ga\in x$ there are coefficients $f_{u,v}\in \F$  for which  $$x\setminus \ga=\sum\,f_{u,v}[u,v]$$ with  summation running over all pairs $(u,v)$ with $|u|+|v|=k-1$ and $\ga,\,\gab\not\in u\cup v.$ 
Then $x=\frac12\big(x\setminus \ga\big)\big((\ga+\gab)+(\ga-\gab)\big)=\frac12\big(\sum\,f_{u,v}[u,v]\big)\big((\ga+\gab)+(\ga-\gab)\big).$ Now observe that 
$$[u,v](\ga+\gab)=[u\cup{\ga}, v] \text{\,\,\,and\,\,\,}[u,v](\ga-\gab)=[u, b\cup{\ga}] $$  if $\ga\in \dot{V}$ and 
$$[u,v](\ga+\gab)=[u\cup{\gab}, v] \text{\,\,\,and\,\,\,}[u,v](\ga-\gab)=[u, b\cup{\gab}]\,$$ if $\ga\in V\setminus \dot{V}.$  \dne

\bigskip 
The incidence relation  for the faces of $\Hy^{n}$ and their pairing give rise to various  incidence maps.
For $x\in L^{n}_{k}$ we set 
\begin{eqnarray} \e(x) &:=& \sum \,y \text{\quad \,with $y\supset x$ and $y\in L^{n}_{k+1},$}\nonumber\\
\p(x) &:=& \sum \,z \text{\quad \,with $z\subset x$ and $z\in L^{n}_{k-1},$}\nonumber\\
\s(x) &:=& \sum_{\gga\in x} \,\big(x\setminus \gga\big) \,\cup\, \ggab\, \nonumber\end{eqnarray}and extend each to a linear map $M^{n}\to M^{n}.$ Evidently the maps restrict to $M_{k}^{n}$ for any $0\leq k\leq n.$ We denote these restrictions by 
$$\e_{k}\!: M^{n}_{k}\to M^{n}_{k+1},\quad \p_{k}\!: M^{n}_{k}\to M^{n}_{k-1}
\text{\,\,\, and \,\,\,}\s_{k}\!: M^{n}_{k}\to M^{n}_{k}.$$

It is easy to check that $\e$ and $\p$ are adjoint to each other in the inner product on $M^{n},$ that is, $(\e(x),y)=(x,\p(y))$ for all $x,\,y\in L^{n}.$ 
Therefore the maps $$\gnp=\p\e\!:\,M^{n}\to M^{n} \text{\quad and \quad} \gnm=\e\p\!:\,M^{n}\to M^{n}$$ are self-adjoint.  Their restrictions are denoted by   $$\gnp_{k}=\p_{k+1}\e_{k}\!:\,M^{n}_{k}\to M^{n}_{k} \text{\quad and \quad} \gnm_{k}=\e_{k-1}\p_{k}\!:\,M^{n}_{k}\to M^{n}_{k}.$$
Similarly,  $\s$ is self-adjoint,  $(\s(x),y)=(x,\s(y))$ for all $x,\,y\in L^{n}.$ Note that if $a,\,b\in L^{n}$ are disjoint sets for which also $a\cdot b$ belongs to $L^{n}$ then we have $$(*)\,\quad \p(a\cdot b)=\p(a)\cdot b+a\cdot\p(b)\text{\quad and \quad}\s(a\cdot b)=\s(a)\cdot b+a\cdot\s(b).$$ 
We will use this fact without further reference. 

\medskip
Our aim is to show that the face modules of $\Hy^{n}$  can be decomposed completely into eigenspaces of these self-adjoint linear maps. Here, by completely we mean that the eigenspaces are in fact irreducible modules for the hyperoctahedral group over any field $\F$ whose  characteristic is $p=0$ or $p>n.$  This will be proved in the next section. It is essential that this decomposition can be computed directly from the incidence relation of the simplicial complex. Our method is entirely general and can be applied to arbitrary ranked partially ordered set. (The resulting decomposition will of course not always yield a decomposition into irreducibles for the automorphism group of the poset.) We start our analysis with the map $\s.$ 

\bigskip
\begin{lem}\label{lem2.2}
  Let $0\leq j\leq k\leq n$ and suppose that $a,\,b$ are disjoint subsets of $\dot{V}$ with  $|a|=k-j$ and $ |b|=j.$ Then $\s_{k}([a,b])=(k-2j)[a,b].$ \end{lem}

\pf For $k=0$ we have  $\s_{0}([\emptyset,\emptyset])=\s_{0}(\emptyset)=0\cdot[\emptyset,\emptyset]$ as required. Suppose therefore that the lemma holds for all $k'<k$ and let $a,\,b$ be disjoint subsets of $\dot{V}$ with $|a|+|b|=k>0$ and $|b|=j.$ 

If $|a|>0$ pick $\ga\in a$ and put $a'=a\setminus \ga.$ Then $[a,b]=(\ga+\gab)\cdot [a',b]$ and therefore $\s([a,b])=\s(\ga+\gab)\cdot [a',b]+(\ga+\gab)\cdot \s([a',b])=(\ga+\gab)\cdot [a',b]+(\ga+\gab)\cdot (k-1-2j)[a',b]$ by induction. Hence $\s([a,b])=[a,b]+(k-1-2j)(\ga+\gab)\cdot [a',b]=(k-2j)[a,b].$ 

If $|b|>0$ pick $\gb\in b$ and put $b'=b\setminus \gb.$ Then $[a,b]=[a,b']\cdot(\gb-\gbb)$ and therefore $\s([a,b])=\s([a,b'])\cdot (\gb-\gbb)+[a,b']\cdot\s(\gb-\gbb)=(k-1-2j+2)[a,b']\cdot(\gb-\gbb)-[a,b']\cdot(\gb-\gbb)$ by induction. Hence  $\s([a,b])=(k-1-2j+2)[a,b]-[a,b]=(k-2j)[a,b].$ \dne

\bigskip
For $0\leq j\leq k$ let  $$C_{k,j}:=\big\{[a,b]\,\,:\,\,a,\,b\,\,\text{are disjoint subsets of $\dot{V}$ with }\,\,\, \,|a|=k-j \text{\,\,and\,\,} |b|=j\,\big\}$$ and let $$E_{k,j}=\F C_{k,j}\,\,\, \text{be the subspace of $M_{k}^{n}$ spanned by $C_{k,j}$}.$$ As before, we set $E_{k,j}=0$ when $C_{k,j}$ is the empty set. 
Together  Proposition~\ref{thm2.1}\, and Lemma~\ref{lem2.2}\, give the following:

\bigskip
\begin{thm}\label{thm2.2}Suppose that $\F$ has characteristic $p\neq 2$ and let $0\leq k\leq n.$  Then \\[-13pt]\begin{eqnarray}\label{decos}
M^{n}_{k}=E_{k,0}\,\oplus \,E_{k,1}\,\oplus\,...\oplus E_{k,k}\,\end{eqnarray}
where $E_{k,j}$ has dimension ${n\choose k}{k\choose j}$  for $0\leq j\leq k$ and basis $C_{k,j}.$
Furthermore, $\s_{k}\!:\,M_{k}^{n}\to M_{k}^{n}$ restricts to $\s_{k,j}\!:\,E_{k,j}\to E_{k,j}$ with $\s_{k,j}=(k-2j) \,I$ for all $0\leq j\leq k$ where $I$ denotes the identity map on $E_{k,j}.$\end{thm}

\medskip
In particular, if $p>n$ then the $\gl_{k,j}=k-2j$ for $0\leq j\leq k$ are the distinct eigenvalues of $\s_{k}$ and \,(\ref{decos})\, is the decomposition of $M^{n}_{k}$ into the eigenspaces. For $p<n$ the eigenspace for $\gl_{k,j}=k-2j$ is of the form $\bigoplus E_{k,j^{*}}$ where the sum runs over all $j^{*}\equiv j
\pmod p.$

\bigskip
Next we consider the interaction between $\p,\,\e$ and $\s.$ The internal structure of a facet (a maximal face) of $\Hy^{n}$ is that of a Boolean algebra on an  $n$-set. Therefore consider also the modules and incidence maps for the Boolean algebra. Recall, we denoted $\dot{V}:=\{\ga_{1},..,\,\ga_{n}\}$ and  $\dot{L}^{n}_{k}:=\{ x\in L^{n}_{k}\,:\, x\subseteq \dot{V}\}.$ Correspondingly we have the facet modules $\dot{M}^{n}_{k}:=\F \dot{L}^{n}_{k}$ and $\dot{M}^{n}:=\oplus \dot{M}^{n}_{k}.$ These are equipped with the same standard inner product as before. The incidence maps are given in the same way:   for $x\in \dot L^{n}_{k}$  we set 
\begin{eqnarray} \de(x) &=& \sum \,y \text{\quad \,with $y\supset x$ and $y\in \dot L^{n}_{k+1},$}\nonumber\\
\dop(x) &=& \sum \,z \text{\quad \,with $z\subset x$ and $z\in \dot L^{n}_{k-1},$}\nonumber \end{eqnarray}
and extend linearly, obtaining the maps $\de,\,\dop\!:\,\dot M^{n}\to \dot M^{n}.$ Note that $\dop$ is the restriction of $\p$ to $\dot M^{n}$ while $\de$ is genuinely different from $\e.$  The restrictions to $\dot M_{k}^{n}$ are denoted by 
$$\de_{k}\!: \dot M^{n}_{k}\to \dot M^{n}_{k+1}\text{\quad and \quad}\dop_{k}\!: \dot M^{n}_{k}\to \dot M^{n}_{k-1} .$$  In the same way we define $$\dgnp=\dop\de\!:\,\dot M^{n}\to \dot M^{n} \text{\quad and \quad} \dgnm=\de\dop\!:\,\dot M^{n}\to \dot M^{n}$$ and for any $0\leq k\leq n$ we have the  restrictions  $$\dgnp_{k}=\dop_{k+1}\de_{k}\!:\,\dot M^{n}_{k}\to \dot M^{n}_{k} \text{\quad and \quad} \dgnm_{k}=\de_{k-1}\dop_{k}\!:\,\dot M^{n}_{k}\to \dot M^{n}_{k}.$$ As before we have  $$\,\quad \dop(a\cdot b)=\dop(a)\cdot b+a\cdot\dop(b)$$
whenever $a$ and $b$ are disjoint subsets of $\dot V.$

\bigskip
The next lemma is crucial for the Boolean algebra.  Suitably formulated it holds in all geometries with Dynkin diagram  of type A.   
 
\medskip
\begin{lem}[A-Type Lemma] \label{atype} Let $0\leq k\leq n.$ Then $$\dgnp_{k}-\dgnm_{k}=(n-2k){I}$$ where ${I}$ denotes the identity map on $\dot M^{n}_{k}.$
\end{lem}

\pf Let $x\in \dot L^{n}_{k}.$ Put $g=\dgnp_{k}-\dgnm_{k}$ and $h=(n-2k){I}.$ It  suffices to show that $(g(x),z)=(h(x),z)$ for every $z\in \dot L^{n}_{k}.$

Since $\dop$ and $\de$ are  adjoint to each other we have    $(g(x),z)=(\de(x),\de(z))-(\dop(x),\dop(z))$ and $(h(x),z)=(n-2k)(x,z).$ First, if $z=x$ then $(g(x),x)=(n-k)-k=(h(x),x)$ as required. Next, if $x\cap z$ has size $k-1$ then $(\de(x),\de(z))=1=(\dop(x),\dop(z))$ and $(\de(x),\de(z))-(\dop(x),\dop(z))=0,$ as required. Finally, if $x\cap z$ has size less than $k-1$ it is easy to check that also here $(g(x),z)=0=(h(x),z).$ \dne 

\medskip
Now we return to the  incidence maps of $\Hy^{n}.$ We observe that these preserve the  eigenspaces of $\s,$ recall the definition of $[a,b]$ earlier: 

\medskip
\begin{lem}\label{lem2.7}
  Let $a$ and $b$ be disjoint subsets of $\dot{V}.$ Then \newline
  (i) \,\, $\p([a,b])=2[\dop(a),b]$ and \newline
  (ii) \,   $\e([a,b])=[\de(a),b].$  
  \end{lem}

\pf By definition we have   $[a,b]=\Pi_{\ga\in a}(\ga+\gab)\cdot\Pi_{\gb\in b}(\gb-\gbb).$ For (i) note that $\p(\ga+\gab)=\emptyset+\emptyset=2\emptyset$ while $\p(\gb-\gbb)=\emptyset-\emptyset=0$ for any $\ga$ and $\gb$ in $\dot V.$ Using \,(*)\, we get $$\p([a,b])=2\sum_{a^{*}}\big(\Pi_{\ga\in a^{*}}(\ga+\gab)\cdot\Pi_{\gb\in b}(\gb-\gbb)\big)=2 \sum_{a^{*}} \,[a^{*},b]$$ where the sum runs over all $a^{*}\subset a$ with  $|a^{*}|=|a|-1.$ Hence $\p[a,b]=2[\dop(a),b].$

For (ii) note that if $x\in L^{n}$ then $\e(x)=x\cdot \sum \gamma$ where the sum runs over all $\gamma$ in $V$ with $\gamma$ not in $x\cup\bar x.$ Therefore $\e(\Pi_{\ga\in a}(\ga+\gab)\cdot\Pi_{\gb\in b}(\gb-\gbb))=\Pi_{\ga\in a}(\ga+\gab)\cdot\Pi_{\gb\in b}(\gb-\gbb)\cdot \sum \gamma$ where the sum runs over all $\gamma$ in $V$ with $\gamma$ not in $a\cup\bar a\cup b\cup \bar b.$ Since by definition $[x,y]=0$ if $x\cap y\neq \emptyset$ we can write $\Pi_{\ga\in a}(\ga+\gab)\,\cdot\,\Pi_{\gb\in b}(\gb-\gbb)\,\cdot\, \sum \gamma= \Pi_{\ga\in a}(\ga+\gab)\,\cdot \sum \gamma\,\cdot\,\Pi_{\gb\in b}(\gb-\gbb)=[\de(a),b].$   \dne 

Let $0\leq j, k\leq n$ and let $E_{k,j}=\F C_{k,j}$ as before, where $E_{k,j}=0$ if $k<j.$ It follows from Lemma~\ref{lem2.7} that the maps $\e_{k}$ and $\p_{k}$ restrict further. We denote these restrictions by   $$\e_{k,j}\!: E_{k,j}\to E_{k+1,j} \text{\,\,\, and \,\,\,} \p_{k+1,j}\!: E_{k+1,j}\to E_{k,j}\,$$ and consequently  we have the maps $$\gnp_{k,j}=\p_{k+1,j}\e_{k,j}\!: E_{k,j}\to E_{k,j} \text{\,\,\, and \,\,\,} \gnm_{k,j}=\e_{k-1,j}\p_{k,j}\!: E_{k,j}\to E_{k,j}.$$

The content of Theorem~\ref{thm2.2} \,and Lemma~\ref{lem2.7}\, can be represented in the following {\sc Figure}~1. It shows the decomposition of $M^{n}_{k}$ for all $0\leq k\leq n$ into submodules and the restrictions of the incidence maps between the corresponding $E_{k,j}.$

{\small
$$\begin{array}{lclclclclclclc}

\,\,\,\,\vdots\!\!\!&&\,\,\,\,\vdots&&\,\,\,\,\vdots& &&&\,\,\,\vdots\,\,&&\,\,\,\vdots&&\,\,\,\,&\vdots\\

M^n_{k+1} \!\!\!&=& E_{k+1,0}&\oplus &E_{k+1,1}&\oplus&\quad\cdots\quad&\oplus &E_{k+1,k}&\oplus &E_{k+1,k+1}\!\!\!&&\!\!\oplus & \, 0\\\\
 {\small \p\!\!\downarrow\!\uparrow\!\!\e}&&{ \ssapp}&&{\ssapp}&&\quad\cdots\quad&&{ \ssapp}& &\ssapp&&&\\\\
M^n_{k} \!\!\!&=&E_{k,0}&\oplus & E_{k,1}&\oplus&\quad\cdots\quad&\oplus & E_{k,k}& 
\oplus &\,\,\,0&&&\\\\
\ssapp&&\ssapp&&\ssapp&&\quad\quad&&\ssapp&&&&&\\
\,\,\,\,\,\vdots\!\!\!&&\,\,\,\,\vdots&&\,\,\,\,\vdots& &\quad\,\,\iddots\quad&&\,\,&&&&&\\\\
M^n_{1} \!\!\!&=&E_{1,0}&\oplus &E_{1,1}&\oplus&\,\,\,\,\,\,\,\,\,\,0&\quad\quad& & &&\\&&&&&&&&&&&&\\
\ssapp&&\ssapp&&\ssapp&&&\quad\quad&  &&\\&&&&&&&&&&&&\\
M^n_{0} \!\!\!&=&E_{0,0}&\oplus &\,\,\,\,0&&&\quad\quad&  &&&&\\
\\
\ssapp&&\ssapp&&&&&\quad\quad&  &&&&\\&&&&&&&&&&&&\\
M^n_{-1} \!\!\!&=&\,\,\,\,0&&&&&\quad\quad&  &&&&\\\\
\end{array}
$$ 

\centerline{\sc Figure 1: Submodules of $M^{n}_{k}$}
}

\bigskip For the hyperoctahedron the next lemma is crucial. It corresponds to the  A-Type Lemma. Suitably formulated it holds in the geometries with Dynkin diagram of type B. 

\medskip
\begin{lem}[B-Type Lemma] \label{btype} Let $0\leq k\leq n.$ Then $$\gnp_{k}-\gnm_{k}=(2n-3k){ I} - \s_{k}$$ where ${I}$ denotes the identity map on $M^{n}_{k}.$
\end{lem} 

\medskip 
\pf Let $x\in L^{n}_{k}.$ Put $g=\gnp_{k}-\gnm_{k}$ and  $h=(2n-3k){I} - \s_{k}.$ It  suffices to show that $(g(x),z)=(h(x),z)$ for every $z\in L^{n}_{k}.$

Since $\p$ and $\e$ are  adjoint to each other we have $(\gnp(x),z)=(\e(x),\e(z))$ and $(\gnm(x),z))=(\p(x),\p(z)).$ Therefore    $$(g(x),z)=(\e(x),\e(z))-(\p(x),\p(z))$$ and $$(h(x),z)=(2n-3k)(x,z)-(\s(x),z).$$ First, if $z=x$ then $(g(x),x)=(2n-2k)-k=(h(x),x)$ as required. Next suppose that $x\cap z$ has size $k-1,$ 
say $x\cap z=\{\gb_{1},..,\gb_{k-1}\}=:b$  so that $x=\chi\cup b$ and  $z=\zeta\cup b$ for some vertices $\chi,\, \zeta\in V.$ Here $(\p(x),\p(z))=1.$ If $\zeta=\bar{\chi}$ then $(\e(x),\e(z))=0$ and so  $(g(x),z)=-1.$ On the right hand side we have   $(h(x),z)=0-(\s(x),z)=-1,$ as required. Otherwise, if $\zeta\neq \bar{\chi}$ then $(\e(x),\e(z))=1$ and so  $(g(x),z)=0.$ On the right hand side we now have   $(h(x),z)=0-(\s(x),z)=0,$ as required. Finally, if $x\cap z$ has size less than $ k-1$ then it is easy to check that $(g(x),z)=0=(h(x),z).$ \dne

Using Lemma~\ref{lem2.2} we note in particular:    

\medskip
\begin{cor}\label{cor2.8} Let $0\leq j\leq k\leq n.$ Then $$\gnp_{k,j}-\gnm_{k,j}=2(n-2k+j)I$$ where ${I}$ denotes the identity map on $E_{k,j}.$
\end{cor}

 The corollary shows the close relationship between the simplex and the hyperoctahedron at the level of the eigenspaces of $\s_{k}\!:$ The A-Type Lemma is a special case of the B-Type Lemma, when $j=0.$  Here the corollary says that $\gnp_{k,0}-\gnm_{k,0}=2(n-2k)I$ which in essence is the statement of the A-Type Lemma, for the following reason. 
The association  $a=\{\ga_{1},..,\ga_{k}\}\longleftrightarrow (\ga_{1}+\bar{\ga}_{1})\cdots(\ga_{k}+\bar{\ga}_{k})=[a,\emptyset]$ \,is a bijection between $\dot{L}^{n}_{k}$ and $C_{k,0}. $ Since  $\dot{L}^{n}_{k}$ and $C_{k,0}$ are bases of $\dot{M}^{n}_{k}$ and $E_{k,0}$  respectively,  we have an $\F{\rm Sym}(\dot V)$-isomorphism $\psi\!:\,\dot{M}^{n}_{k}\to E_{k,0}.$   The factor $2$ in the corollary is accounted for by Lemma~\ref{lem2.7}(i). This isomorphism is not obvious from the face complex of $\Hy^{n}.$

\bigskip
The corollary can be used  to decompose  the modules in each column of {\sc Figure}~1\, further. Specifically,  we are interested in the eigenspaces of $\gnp_{k,j}\!:\,E_{k,j}\to E_{k,j}$ for each $0\leq j\leq k.$ We need a general fact about eigenvalues and eigenspaces of  the product of two linear maps.  The linear map $\chi:U\to U$ of the vector space $U$ is said to be {\it semi-simple} if $U$ has an eigenbasis for $\chi.$ Equivalently, $\chi$ is semi-simple  if and only if the minimum polynomial of $\chi$ is a product of distinct linear terms.

\medskip
\begin{lem}\label{lem2.9} Let $U$ and $V$ be finite-dimensional vector spaces and suppose that $\varphi :U\to V$ and $\psi  : V\to U$ are linear maps. Then the following hold:\\\\[-6pt](i) \,\,\, The maps $\psi \circ\varphi :U\to U$ and $\varphi \circ\psi  :V\to V$ have the same non-zero eigenvalues. If $\gl$ is a nonzero eigenvalue  let  $U_{\gl}\subseteq U$  and $V_{\gl}\subseteq V$ denote the eigenspaces of \,$\psi \circ\varphi$\, and \,$\varphi \circ\psi,$\, respectively. Then $\varphi $ and $\psi  $ restrict to  isomorphisms $\varphi :U_{\gl}\to V_{\gl}$ and $\psi  : V_{\gl}\to U_{\gl}.$  \\\\ [-8pt]
(ii) \,\,\,Suppose that  $\psi \circ\varphi\!:\, U\to U$ is semi-simple. Then   $\varphi \circ\psi\!:\, V\to V$ is semi-simple   unless $0$ is an eigenvalue of  both maps and $\varphi \circ\psi$ has minimum polynomial $(x-\gl_{1})\cdots (x-\gl_{t})\cdot x^{2}$ where $\gl_{1},\dots,\,\gl_{t}$ are the non-zero eigenvalues of $\psi \circ\varphi$ and $\varphi \circ\psi.$ \end{lem}

The exception does occur. It is  easy to find $2\times 2$ matrices $A$ and $B$ for which$AB=0$ while $BA\neq 0$ and $(BA)^{2}=0.$

\medskip 
\pf (i) Let $f$ be an eigenvector of $\psi \circ\varphi $ for the eigenvalue $\gl\neq 0.$  Then $(\psi\circ\varphi)(f)=\gl\,f$ and so $(\varphi\circ\psi  )(\varphi (f))=\lambda\varphi (f)$ where $\varphi (f)\neq 0$ since $\gl\neq 0.$ Hence $\gl$ is an eigenvalue of  $\varphi\circ\psi.$ Furthermore, the map $\varphi :U_{\gl}\to V_{\gl}$ is injective since $\gl\neq 0.$ Now  reverse the order of $\psi  $ and $\varphi .$ 

(ii) It is well-known that $\psi \circ\varphi $ and $\varphi\circ \psi $ have the same characteristic polynomial, apart from a factor $x^{m}$ with $m=\pm (\dim U-\dim V).$ By assumption therefore the minimum polynomial of  $\varphi\circ \psi $ is a product of linear terms. Let $\gl$ be an eigenvalue of $\varphi \circ\psi$ and suppose that there is some $v\in V$ such that $(\varphi \circ\psi-\gl)^{2}(v)=0$ but $(\varphi \circ\psi-\gl)(v)\neq 0.$

Expand $\psi(v)\in U$ into eigenvectors of $\psi \circ\varphi.$ Since $\psi\circ(\varphi \circ\psi-\gl)^{2}(v)=(\varphi \circ\psi-\gl)^{2}(\psi(v))=0$ it follows that $\psi(v)=0$ or $\psi(v)\neq 0$ is an eigenvector of $\psi\circ \varphi$ for the eigenvalue $\gl.$ In either case, $0=(\varphi \circ\psi-\gl)^{2}(v)=(\varphi \circ(\psi\circ\varphi)\circ\psi)(v)-2\gl(\varphi \circ\psi)(v)+\gl^{2}(v)=-\gl(\varphi \circ\psi)(v)+\gl^{2}(v)=-\gl(\varphi \circ\psi-\gl)(v).$ Therefore $\gl=0,$ since $(\varphi \circ\psi-\gl)(v)\neq 0.$ In particular, $(x-\gl_{i})^{2}$ does not divide the minimum polynomial of $\varphi\circ \psi $ unless $\gl_{i}=0.$

If $\gl=0$ is not an eigenvalue of $\psi\circ \varphi$ then $\psi(v)=0,$ contradicting $(\varphi \circ\psi-\gl)(v)\neq 0.$ In particular, $x^{2}=(x-\gl)^{2}$ does not divide the minimum polynomial of $\varphi\circ \psi.$ Finally suppose that $(\psi\circ \varphi)^{s}(v)=0$ for $s\geq 2.$ Then $0=(\psi\circ \varphi)^{s-2}\circ\psi\circ(\varphi\circ\psi)(\varphi(v))$ and since $(\varphi\circ\psi)(\varphi(v))=0$ we have that $x^{3}$ does not divide the minimum polynomial of $\varphi\circ \psi.$  \dne


In order to state the main result of this section we recall some earlier details. For a field $\F$ of characteristic $p\neq 2$ and integers $0\leq j\leq k\leq n$ we denoted by  $E_{k,j}=\F C_{k,j}\subseteq M^{n}_{k}$ the vector space with basis $C_{k,j}.$ Then  \\[-13pt]\begin{eqnarray}\label{decos}
M^{n}_{k}=E_{k,0}\,\oplus \,E_{k,1}\,\oplus\,...\oplus E_{k,k}\,\end{eqnarray}
where $E_{k,j}$ has  dimension \\[-15pt] $$d_{k,j}:=\dim E_{k,j}={n\choose k}{k\choose j}.$$\\[-5pt]  For convenience we put $d_{k,j}=0 $ if $k<j$ or $n<k.$ 
By Theorem~\ref{thm2.2}\,  $\s_{k}\!:\,M^{n}_{k}\to M^{n}_{k}$  has eigenvalue $k-2j$ on $E_{k,j}.$

\bigskip
\begin{thm}[Spectral Decomposition for ${\mathcal H}^{n}$]\label{fulldecomp} Let $\F$ be a field of characteristic $p>n$ or $p=0.$ Suppose that  $2\leq n$ and  $0\leq k\leq n.$  Then the following hold for every  $0\leq j\leq k$ and $m:=\min\{k-j, \,n-k\}.$:\\[10 pt]
(i) \,\,\, The eigenvalues of  $\gnp_{k,j}\!:\,E_{k,j}\to E_{k,j}$ are $$\gl_{k,j,i}:=2(k-j+1-i)(n-k-i)$$ with $0\leq i\leq m.$ In particular, the $\gl_{k,j,i}$ are pairwise distinct for $0\leq i\leq m$ with $\gl_{k,j,i}=0$ if and only if $i=n-k\leq k-j.$ \\[10 pt]
 (ii)\,\,  The map $\gnp_{k,j}\!:\,E_{k,j}\to E_{k,j}$ is semi-simple. Let $0\leq i\leq m$ and let $E_{k,j,i}\subseteq E_{k,j}$ denote the eigenspace for $\gl_{k,j,i}.$ Then $\dim(E_{k,j,i})=d_{j+i,j}-d_{j+i-1,j}.$  \\[10 pt] 
 (iii) \, For every  $0\leq i\leq n-j$ the maps $\e$ and $\p$  restrict to isomorphisms $$\e, \p\!:\,E_{j+i,j,i}\,\,\longleftrightarrow\,\, E_{j+i+1,j,i}\,\,\longleftrightarrow\,\, \dots \,\,\longleftrightarrow\,\,  E_{n-i,j,i}$$  while $\p(E_{j+i,j,i})\,=\,0\,=\,\e(E_{n-i,j,i}).$ \end{thm}
 
 \bigskip
 \pf (i) Fix $0\leq j\leq n$ and proceed  by induction over $k,$ with $j\leq k\leq n.$ For $k=j$ we have  $i=0$ and  $\gnm_{j,j}\!:\,E_{j,j}\to E_{j,j}$ is the zero map since $E_{j-1,j}=0$ by definition. By Corollary~\ref{cor2.8}  therefore $\gnp_{j,j}\!:\,E_{j,j}\to E_{j,j}$ is constant,  $\gnp_{j,j}=2(n-j)I=2(k-j+1-i)(n-k-i)I$ for $i=0,$  as claimed. 

Now assume that the statement in (i) holds for all values less than $k$ and let $\gl$ be an eigenvalue of  $\gnp_{k,j}\!:\,E_{k,j}\to E_{k,j}.$ Then $\gl'=\gl-2(n-2k+j)$ is an eigenvalue of $\gnm_{k,j}\!:\,E_{k,j}\to E_{k,j}$ by Corollary~\ref{cor2.8}. Now distinguish two cases, according to Lemma~\ref{lem2.9}.

If $\gl'$ is an eigenvalue of  $\gnp_{k-1,j}\!:\,E_{k-1,j}\to E_{k-1,j}$ then by induction $$\gl'=2((k-1)-j+1-i)(n-(k-1)-i)$$ for some $0\leq i\leq m'=\min\{(k-1)-j, \,n-(k-1)\}.$ Notice that \begin{eqnarray}\label{ind2}
&&[2((k-1)-j+1-i)(n-(k-1)-i)]+[2(n-2k+j)]\nonumber\\
&&=2(k-j+1-i)(n-k-i)\end{eqnarray} for all $i,\,j,\,k,\,n.$ Hence the required property holds, in particular that $0\leq i\leq \min\{k-j,\,n-k\},$ unless  $i=n-(k-1)\leq (k-1)-j.$ But in this case  $\gl'=0$ and $\gl=2(n-2k+j)=2(k-j+1-i)(n-k-i)$ for $i=k-j.$ Note that $i\leq n-k$ since $(k-1)-j\geq n-(k-1).$ This concludes the proof of (i) when $\gl'$ is an eigenvalue of  $\gnp_{k-1,j}.$

Otherwise $\gl'$ is not an eigenvalue of  $\gnp_{k-1,j}\!:\,E_{k-1,j}\to E_{k-1,j}$ and hence $\gl'=0$ by Lemma~\ref{lem2.9}. In particular, $\e_{k-1}\!:\,E_{k-1,j}\to E_{k,j}$ is injective as otherwise $0$ is an eigenvalue of $\p_{k}\e_{k-1}.$ Therefore $\dim(E_{k-1,j})\leq \dim(E_{k,j})$ and from Theorem~\ref{thm2.2}\, we have $k-j\leq n-k+1.$ 

Since $\gl'=0$ we have $\gl=2(n-2k+j)=2(k-j+1-i)(n-k-i)$ for $i=k-j$ as above.  This completes the proof unless $n-k<k-j.$ In this remaining case we have $k-j=n-k+1$ which implies that $\gl=2(n-2k+j)=-2\neq 0$ and $\dim(E_{k-1,j})=\dim(E_{k,j})$ by Theorem~\ref{thm2.2}. This case does not occur: If $\e_{k-1}\p_{k}(f)=0$ for some $0\neq f\in E_{k,j}$ then $\p_{k}(f)=0$ since $\e_{k-1}$ is injective. Since $\e_{k+1}$ is also surjectice there is some $f'\in E_{k-1,j}$ with $f=\e_{k-1}(f')$ and so $\p_{k,j}\e_{k-1,j}(f')=0$ while $f'\neq 0,$ contradicting the assumption that $\gl'=0$ is not an eigenvalue of  $\gnp_{k-1,j}\!:\,E_{k-1,j}\to E_{k-1,j}.$  The remaining assertions in (i) are immediate.

\medskip
(ii) We show that $\gnp_{k,j}\!:\,E_{k,j}\to E_{k,j}$ is semi-simple for all $0\leq j\leq k\leq n.$  As before fix $0\leq j\leq n$ and proceed  by induction over $k.$ For $k=j$ we have  $i=0$ and  $\gnm_{j,j}\!:\,E_{j,j}\to E_{j,j}$ is the zero map, hence  semi-simple.  The dimension of $E_{j,j}=E_{j,j,0}$ is $d_{j,j,0}=d_{j,j}-0=d_{j,j}-d_{j-1,j}.$

Next assume that the result holds for all values $<k.$ Let $\gl_{k-1,j,i}$ and $\gl_{k,j,i}$ be the eigenvalues from (i) and denote the corresponding eigenspaces by $E_{k-1,j,i}$ and $E_{k,j,i}$ respectively. By induction $E_{k-1,j}=E_{k-1,j,0}+...+E_{k-1,j,m'}$ with $m'=\min\{(k-1)-j,\,n-(k-1)\}.$ Since $\gnm_{k,j}=\gnp_{k,j}-2(n-2k+j)$ by the B-Lemma it is sufficient to show that $\gnm_{k,j}\!:\,E_{k,j}\to E_{k,j}$ is semi-simple and that $E_{k,j,i}$ has dimension $d_{k,j,i}$ for all  $0\leq i\leq m.$

Recall that $\gnm_{k,j}=\e_{k-1,j}\p_{k,j}$  and  $\gnp_{k-1,j}=\p_{k,j}\e_{k-1,j}$ with  $\e_{k-1}:E_{k-1,j}\to E_{k,j}$ and $\p_{k}:E_{k,j}\to E_{k-1,j}.$ We can now apply Lemma~\ref{lem2.9}.

First suppose that $0$ is not an eigenvalue of $\gnp_{k-1,j}.$  By the first part of the theorem  this happens if and only if $n-(k-1)>(k-1)-j.$  In particular, $\e_{k-1}:E_{k-1,j}\to E_{k,j}$ is injective, $E_{k,j,i}=\e_{k-1,j}(E_{k-1,j,i})$ is an isomorphism by Lemma~\ref{lem2.9}\,
and $\dim(E_{k,j,i})=\dim(E_{k-1,j,i})=d_{j+i,j}-d_{j+i-1}$ for $i<m$. For $i=m=k-j$ we have $E_{k,j,i}={\rm ker}(\p_{k,j})$ and $d_{k,j,i}=d_{k,j}-d_{k-1,j}.$ As required, $\gnm_{k,j}\!:\,E_{k,j}\to E_{k,j}$ is semi-simple and $E_{k,j,i}$ has dimension $d_{k,j,i}$ for all  $0\leq i\leq m.$

Secondly suppose that $0$ is an eigenvalue of $\gnp_{k-1,j},$  equivalently that  $n-(k-1)\leq (k-1)-j.$ If $0$ is not an eigenvalue of $\gnm_{k,j}$ then $\gnm_{k,j}\!:\,E_{k,j}\to E_{k,j}$ is semi-simple by Lemma~\ref{lem2.9}(ii) and it is straightforward to show that all eigenspaces have the given dimension. It remains to consider the case when $0$ is an eigenvalue of $\gnm_{k,j}$ when equivalently $2(n-2k+j)$ is an eigenvalue of $\gnp_{k,j.}$  The equation $2(n-2k+j)=2(k-j+1-i)(n-k-i)$ has two solutions, $i=k-j$ and $i=n-k+1.$ But both contradict the requirement $i\leq \min\{n-k,k-j\}$ and $n-(k-1)\leq (k-1)-j.$ Hence $0$ is not an eigenvalue of both $\gnm_{k,j}$ and $\gnp_{k-1,j}.$

\medskip
(iii) This part follows from Lemma~\ref{lem2.9} and the first part of the theorem. Observe that  that $\gl_{j+i,j,i}\,,\,\gl_{j+i+1,j,i}\,,...,\,\,\gl_{n-i-,j,i}$ are all non-zero while $\p(E_{j+i,j,i})\,=\,0\,=\,\e(E_{n-i,j,i})$ by definition. \dne
 

\bigskip
{\sc Comments:\,}  (1) From the definition, if  $\gnp_{k,j}\!:\,E_{k,j}\to E_{k,j}$ is semi-simple then 
$$E_{k,j}=E_{k,j,1}\,\oplus\,E_{k,j,2}\,\oplus\,...\oplus E_{k,j,m}$$ 
is the  decomposition into the eigenspaces of $\gnp_{k,j}.$ From this we have the decomposition $$M^{n}_{k}\,=\,\,\bigoplus_{j} \,E_{k,j}\,=\,\,\bigoplus_{j,i} \,E_{k,j,i}\quad\text{and}\quad M^{n}\,=\,\,\bigoplus_{k,j,i} \,E_{k,j,i}
$$ into eigenspaces of $\s_{k}$ and  $\gnp_{k}.$  Each element $f$ in  $M_{k}^{n}$  can be decomposed as $$f=\sum_{j=0}^{k}\,f_{k,j}\text{\quad and \quad}f_{k,j}=\sum_{i=0}^{m}\,f_{k,j,i}\label{deco}$$ with uniquely determined $f_{k,j}\in E_{k,j}$ and $f_{k,j,i}\in E_{k,j,i}.$ We call $f_{k,j}$ and $f_{k,j,i}$ the {\it spectral components}, or just {\it components} of $f. $ In the next section we sketch  how these components are computed from the eigenvalues of $\s$ and $\gnp.$ 

(2) The decomposition can be illustrated by selecting a column in  {\sc Figure}~1\, and expanding each term further into eigenspaces of $\gnp$ with isomorphisms between them as provided by $\e$ and $\p.$ This is shown in {\sc Figure}\, 2.  Informally,  each columns is symmetrical about its middle position at $\frac{n+j}{2}.$ In Section~4 we show that  $\e,\,\p\!:\, M_{k}^{n}\leftrightarrow M^{n}_{k+1}$ are in some sense the only maps that commute with the group of all  automorphisms  of $\Hy^{n}.$   \\[-4pt]

 {\small
$$\begin{array}{lclclclclclclc}
E_{n+1,j} \!\!\!&=& \,\,0& &&&&&&&&&&\\\\
{ \ssapp}&&{ \ssapp}&&\,\,\,&&&&& &&&&\\\\
E_{n,j} \!\!\!&=& E_{n,j,0}&\oplus &\quad 0&&&&&&&&&\\\\
{ \ssapp}&&{ \ssapp}&&\,\,\,{ \ssapp}&&&&& &&&&\\\\
E_{n-1,j} \!\!\!&=& E_{n-1,j,0}&\oplus &E_{n-1,j,1}&\oplus&\,\,\,0\quad\quad&&&&&&&\\\\{ \ssapp}&&{ \ssapp}&&{\ssapp}&&{\ssapp}&&& &&&&\\
\,\,\,\,\,\vdots\!\!\!&&\,\,\,\,\vdots&&\,\,\,\vdots& &\,\,\,\,\vdots&&\ddots\quad&&&\\\\\\\\
E_{j+k,j} \!\!\!&=&E_{j+k,j,0}&\oplus &E_{j+k,j,1}&\oplus&\,\,\,\cdots &\oplus &E_{j+k,j,k}  &\,\oplus&\,\,\,0\\&&&&&&&&&&\\
{ \ssapp}&&\ssapp&&\ssapp&&&& \,\,\ssapp&&\\
\,\,\,\,\,\vdots\!\!\!&&\,\,\,\,\vdots&&\,\,\,\,\vdots& &\,\,\iddots&&&&&\\\\
E_{j+1,j} \!\!\!&=&E_{j+1,j,0}&\oplus &E_{j+1,j,1}&\oplus&\,\,\,\,\,0\,&& & &&\\&&&&&&&&&&\\
{\ssapp}&&\ssapp&&\,\,\ssapp&&&&  &&\\&&&&&&&&&&\\
E_{j,j} \!\!\!&=&E_{j,j,0}& \oplus&\,\,\,\,0\,\,&&&&  &&\\\\
{\ssapp}&&\ssapp&&&&&&  &&\\&&&&&&&&&&\\
E_{j-1,j} \!\!\!&=&\,\,0&&\,\,\,\,\,&&&&  &&\\
\\
\end{array}
$$ \\[-3pt]
\centerline{\sc Figure 2: Submodules of $E_{k,j}$}
}

\bigskip
(3) We have  mentioned earlier (in the paragraph following  Corollary~\ref{cor2.8}) that  the column $j=0$ in {\sc Figure}~1 corresponds to the modules of the Boolean lattice $\dot{L}^{n}.$  Therefore Theorem~\ref{fulldecomp} and the  isomorphism $\psi\!:\,\dot{M}^{n}_{k}\to E_{k,0}$ discussed before give the decomposition for the modules of the Boolean algebra $$\dot{M_{k}}=\dot{E}_{k,0}\oplus \dot{E}_{k,1}\oplus \cdots \oplus \dot{E}_{k,k}$$
where  $\dot{E}_{k,i}=\psi^{-1}(E_{k,0,i})$ for $0\leq i\leq k.$  The same decomposition can be obtained directly from the A-Type Lemma by the method of the proof above. For reference we state this special case separately: 

\bigskip
\begin{thm}[Spectral Decomposition for Boolean Algebra]\label{Bfulldecomp}  Let $\F$ be a field of characteristic $p>n$ or $p=0.$ Suppose that $2\leq n$ and $0\leq k\leq n,$ and put $m:=\min\{k, \,n-k\}.$  Then the following hold:\\[10 pt]
(i) \,\,\, The eigenvalues of  $\dgnp_{k}\!:\,\dot M_{k}\to \dot M_{k}$ are $$\gl_{k,i}:=(k+1-i)(n-k-i)$$ with  $0\leq i\leq m.$ In particular, the $\gl_{k,i}$ are pairwise distinct for $0\leq i\leq m,$ with $\gl_{k,i}=0$ if and only if $i=n-k\leq k.$  \\[10 pt]
 (ii)\,\,  The map $\dgnp_{k}\!:\,\dot M_{k}\to \dot M_{k}$ is semi-simple. Let $0\leq i\leq m$ and let $\dot E_{k,i}\subseteq \dot M_{k}$ denote the eigenspace for $\gl_{k,i}.$ Then $\dim(\dot E_{k,i})={n\choose i}-{n\choose i-1}.$  \\[10 pt] (iii) \, For every  $0\leq i\leq n$ the maps $\e$ and $\p$  restrict to isomorphisms $$\e, \p\!:\,\dot E_{i,i}\,\,\longleftrightarrow\,\, \dot E_{i+1,i}\,\,\longleftrightarrow\,\, \dots \,\,\longleftrightarrow\,\, \dot E_{n-i,i}$$  while $\p(\dot E_{i,i})\,=\,0\,=\,\e(\dot E_{n-i,i}).$ \end{thm}

\medskip 
 The proof is immediate from Theorem~\ref{fulldecomp} and Lemma~\ref{lem2.7}, alternatively adapt the proof of Theorem~\ref{fulldecomp} using the A-Type Lemma directly.
 
 \medskip
(3) Similar arguments can be used to obtain such decompositions when the characteristic $p$  of $\F$ satisfies  $2<p\leq n.$ The situation is more involved due to the fact that  eigenvalues coincide for distinct parameter triples $k,j,i.$

\bigskip
\section{\sc Applications
}

The decomposition of the face modules of the hyperoctahedron has many applications, we mention some of these now. The first is a rank formula for the incidence matrices ${\mathcal S}_{t,k}$ between the $t$- and $k$-faces of $\Hy^{n}.$ Next we determine the projection maps $M_{k}\to E_{k,j}$ and $E_{k,j}\to E_{k,j,i}.$ We show that these maps can be computed explicitly from the eigenvalues and so the spectral decomposition of elements in $M_{k}$ can be written down. For instance, if $B$ is a set of $k$-faces we can consider the element 
$$f_{B}\,:=\,\sum_{x\in B}\,x\,=\sum_{j}\,f_{k,j}\,=\,\sum_{j,i}\,f_{k,j,i}\quad\text{in} \,\,\,M^{n}_{k}$$ and its spectral components. It turns out that quite basic  information about the components of $f_{B}$ -- such as their euclidean norm -- determines the `shape' of  $B.$  The notion of shape becomes more evident when we discuss $t$-designs on the hyperoctahedron and orbit theorems for groups of automorphisms of $\Hy^{n}.$ In the same context we mention also the connection to association schemes on the hyperoctahedron. In this section all modules and maps refer to $\Hy^{n}$ and are defined over a coefficient field $\F$ of characteristic $p>n$ or $p=0.$

\bigskip\medskip

{\sc 3.1 Incidence Rank:} \, Let $0\leq t < k\leq n$ and suppose that $\F$ is a field of characteristic $p=0$ or $p>n.$ Put 
$$\e_{t}^{k}\!:=\,\frac{1}{(k-t)!}\,\,\,\e_{k-1}\circ \e_{k-2}\circ ... \circ \e_{t}\!:\,M_{t}\to M_{k}.$$ Then $\e_{t}^{k}(x)=\sum y$ where the sum runs over all elements $y\in L^{n}_{k}$ so that $x\subseteq y.$ (The term $(k-t)!$ counts the number of distinct chains  $x=x_{t}\subset x_{t+1}\subset ... \subset x_{k}=y$ with $x_{\ell}\in L^{n}_{\ell}.)$ The map $\e^{k}_{t}$ is the incidence map of the incidence structure $(L^{n}_{t},\,L^{n}_{k},\subseteq). $

\bigskip

\begin{thm}\label{rank} Suppose that  $0\leq j\leq t\leq k\leq n$ and that  $\F$ is a field of characteristic $p=0$ or $p>n.$ Then $\e_{t}^{k}$ restricts to a map $E_{t,j}\to E_{k,j}$ and this restriction has rank $$\min\Big\{{n\choose t}{t\choose j},\,{n\choose k}{k\choose j}\Big\}$$ over $\F.$  In particular,  $\e_{t}^{k}\!:\,M_{t}\to M_{k}$ has rank $$\sum_{0\leq j\leq t}\min\Big\{{n\choose t}{t\choose j},\,{n\choose k}{k\choose j}\Big\}$$  over $\F.$ 
\end{thm}

{\small \sc Comments:} 1. The matrix of $\e_{t}^{k}\!:\,M_{t}\to M_{k}$ with respect to the bases $L_{t}^{n}$ and $L^{n}_{k}$ is  the $(t,k)$-incidence matrix ${\mathcal S}_{t,k}$ which we defined in Section $1.$ Its rows and columns are indexed by the $t$- and $k$-faces of ${\mathcal H}^{n},$  respectively, where the $(x,y)$-entry of ${\mathcal S}_{t,k}$ is $1$ if $x\subseteq y$ and $0$ otherwise.   The second part of the theorem is therefore   Theorem~1.1 \,in the Introduction.

2. The first part of the theorem can be stated as saying that  $\e_{t}^{k}\!:\,E_{t,j}\to E_{k,j}$ has maximum rank, equal to $\min\{\dim(E_{t,j}),\,\dim(E_{k,j})\}.$
For instance, according to our earlier comments $\e_{t}^{k}\!:\,E_{t,0}\to E_{k,0}$ is essentially the incidence map between the $t$- and $k$-subsets of an $n$-set. Therefore we have the well-known fact that in the Boolean lattice  the incidence maps have maximum rank.  Rank maximality is common in many incidence structures, including finite projective spaces, see~\cite{MS2}.

However, this property does not transfer to the incidence maps between  $t$- and $k$-faces of the octahedron in general. For instance, $\Hy^{3}$ has $12$ edges and $8$ faces while  ${\mathcal S}_{2,3}$ has only rank~$7.$ The theorem however does show that  
$${\rm rank}({\mathcal S}_{t,k})=\min\{\dim(M^{n}_{t}),\,\dim(M^{n}_{k})\}$$ when  $t+k\leq n,$ and this is not difficult to verify from the theorem.  In fact, this property holds for the corresponding incidence matrices in any pure simplicial complex, see~\cite{MS2}.

\medskip
{\it  Proof of Theorem~\ref{rank}:}\, Recall that $\dim(E_{s,j})={n\choose s }{s\choose j}$ for $0\leq j\leq s\leq n.$
Let $r$ be the rank of $\e_{t}^{k}\!:\,E_{t,j}\to E_{k,j}.$ Then   $r\leq  \min\{{n\choose t}{t\choose j},\,{n\choose k}{k\choose j}\}.$ First suppose that ${n\choose t}{t\choose j}\leq {n\choose k}{k\choose j},$ equivalently $(n-k)!(k-j)!\leq (n-t)!(t-j)!$. It follows that if $0\leq i\leq \min\{n-t, t-j\}$ then $0\leq i\leq \min\{n-k, k-j\},$ since $0\leq j\leq t\leq k.$ Therefore by Theorem~\ref{fulldecomp}, if $0\neq E_{t,j,i}$, then $\e_{t}^{k}:E_{t,j,i}\to E_{k,j,i}$ is an isomorphism. Hence, for every summand in $E_{t,j}=\bigoplus_{i} E_{t,j,i}$ an isomorphic copy appears in  $\e_{t}^{k}(E_{t,j})\subseteq E_{k,j}.$ In particular,  $r=\dim(E_{t,j}).$ The same argument works also in the second case, when ${n\choose t}{t\choose j}> {n\choose k}{k\choose j}.$ The second part of the theorem is evident since $M^{n}_{t}=\bigoplus_{j=0}^{t}\,E_{t,j}.$ \dne

\bigskip
{ \sc 3.2 Projection onto Eigenspaces:}  We turn to the computation of the decompositions $M^{n}_{k}=E_{k,0}\oplus...\oplus E_{k,k}$ and $E_{k,j}=E_{k,j,0}\oplus...\oplus E_{k,j,m}.$

Let $\psi\!:\,U\to U$ be a semi-simple linear map with distinct eigenvalues $\gl_{0},..,\,\gl_{m}$ and eigenspaces $U_{0},..,\,U_{m}.$ Hence $U=U_{0}\oplus ..\oplus U_{m}.$ Let $\pi_{i}\!:\,U\to U_{i}$ for $0\leq i\leq m$ be the corresponding projection (minimal idempotent). Thus $\pi_{i}^{2}=\pi_{i},$ $\pi_{i}\pi_{j}=0$ for $i\neq j$ and $\sum_{i=1}^{m}\,\pi_{i}={I},$ the identity map $U\to U.$ Let $$\mu_{\psi}(x):=(x-\gl_{0})\cdots(x-\gl_{m})$$ be the minimum polynomial of $\psi$ and for $0\leq j\leq m$ let $$\mu_{\psi,\,j}(x):=(x-\gl_{0})\cdots (x-\gl_{j-1})\cdot (x-\gl_{j+1})\cdots (x-\gl_{m})$$ be the  polynomial with  $(x-\gl_{j})\mu_{\psi,\,j}(x)=\mu_{\psi}(x). $ The following is an essential part of the spectral theorem in linear algebra which is commonly not stated in the literature: 

\medskip
\begin{lem}\label{proj} Let $\psi\!:\,U\to U$ be semi-simple with minimum idempotents $\pi_{j}\!:\,U\to U_{j}$ for $0\leq j\leq m.$ Then $\pi_{j}=[\mu_{\psi,\,j}(\gl_{j})]^{-1}\,\mu_{\psi,\,j}(\psi).$
\end{lem}

\pf  Let $\pi_{i}=[\mu_{\psi,\,i}(\gl_{i})]^{-1}\,\mu_{\psi,\,i}(\psi).$ Then we have  $\mu_{\psi,j}(\psi)\pi_{i}=0$ for $i\neq j$ and $\mu_{\psi,j}(\psi)\pi_{j}=\mu_{\psi,j}(\gl_{j})\pi_{j}.$ Hence if $f=\pi_{1}(f)+\dots +\pi_{m}(f)$ is an element in $U$ then $$\mu_{\psi,j}(f)=\mu_{\psi,j}(\sum_{i}\, \pi_{i}(f))=\mu_{\psi,j}( \pi_{j}(f))=\mu_{\psi,j}(\gl_{j})\pi_{j}(f)$$ as required. \dne

From this lemma and the eigenvalues $\gl_{k,j}$ of $\s$ given in Theorem~\ref{thm2.2} we compute the projections $\pi_{k,j}\!:M_{k}\to  E_{k,j}$ for $0\leq j\leq k$ as 
\begin{eqnarray}\label{dec111}\pi_{k,j}=\prod_{0\leq s\leq k,\,s\neq j}\,\frac{1}{2(s-j)}\,(\s_{k}-(k-2s)I).
\end{eqnarray}
where $I$ denotes the identity map on the corresponding space. Similarly, for fixed $0\leq j\leq k\leq n$ let $m:=\min\{k-j,\,n-k\}.$ From  the eigenvalues $\gl_{k,j,i}$ of $\gnp_{k,j}$ in Theorem~\ref{fulldecomp} we have the projections $\pi_{k,j,i}\!:E_{k,j}\to  E_{k,j,i}$ for $0\leq i\leq m$ given by
\begin{eqnarray}\label{dec112}\pi_{k,j,i}=\prod_{0\leq t\leq m,\,t\neq i}\,\frac{1}{2(t-i)(n-j+1-2i-t)}\,(\gnp_{k,i}-2(k-j-t+1)(n-k-t)I).
\end{eqnarray}
Therefore we have

\bigskip
\begin{prop}\label{compo} Let $0\leq k\leq n$ and let $f$ be an element in $M_{k}^{n}.$ Then the spectral components $f=f_{k,0}+...+f_{k,k} $ (with $f_{k,j}\in E_{k,j}\,\,\,\text{for}\,\,\, 0\leq j\leq k)$ and $f_{k,j}=f_{k,j,0}+...+f_{k,j,m} $ (with $f_{k,j,i}\in E_{k,j,i}\,\,\,\text{for}\,\,\, 0\leq j\leq m:=\min\{k-j,\,n-k\})$
are given by $$f_{k,j}=\pi_{k,j}(f) \quad\text{and}\quad  f_{k,j,i}=\pi_{k,j,i}(f_{k,j})$$ where  $\pi_{k,j}$ and $\pi_{k,j,i}$ are as  in \,(\ref{dec111})\, and \,(\ref{dec112}).
\end{prop}

\bigskip

\bigskip
In some situations the components of an element in $M^{n}_{k}$ play an important role combinatorially. In particular, if $B$ is a family of $k$-faces of $\Hy^{n}$ then we may consider the element $f=f_{B}$ given by   
   $$f\,:=\,\,\sum_{x\in B}\, x\,\,=\,\,\sum_{j,i}\,f_{k,j,i}$$ with spectral components $f_{k,j,i}.$ If we denote $l_{k}:=\sum \,x$ (with $x\in L_{k}^{n})$ then $(f,\,l_{k})=|B|$ and the  orthogonality of eigenvectors for different eigenvalues gives $(f\,,\,l_{k}) =(f_{k,0,0}\,,\,l_{k}).$ Since $E_{k,0,0}$ is one-dimensional $f_{k,0,0}$ is a multiple of $l_{k}$ and since  $(l_{k}\,,\,l_{k})=|L^{n}_{k}|$ we have 
 \begin{eqnarray}\label{deco02}f_{k,0,0}=\frac{|B|}{|L^{n}_{k}|}\,\,l_{k}\,.\,\end{eqnarray}
Therefore $f_{k,0,0}$ and $|B|$ determine one another, for arbitrary $B.$  Similarly,  
\begin{eqnarray}\label{deco03}|B|=(f,f)=\sum_{j,i}\,||f_{k,j,i}||^{2}\,\end{eqnarray} from the orthogonality of spectral components. Some of the other spectral components of $f_{B}$ carry more specific combinatorial information about $B,$ as the following shows.

\bigskip
{\sc 3. Designs:} Let $V$ be a set of $v$ points and let $0<t\leq k\leq v.$ Suppose that $B$ is a family of $k$-element subsets of $V.$ Then $(V,\,B)$ is a  $t-(v,k,\ell)$-design if every  $t$-element subset of $P$ is contained in exactly $\ell$ members of $B.$ This is the usual definition of a block design. Several variations are considered in the literature, including designs on finite projective or affine spaces. The following is the natural extension to partially ordered sets generally: 

{\small \sc Definition:}  Let $({\mathcal L},\leq )$ be a  ranked partially ordered set. Then a $(t,k,\ell)$-design on ${\mathcal L}$ is a set $B\subseteq L_{k}$ so that for every $x\in L_{t}$ there are exactly $\ell$ members $y$ of $B$ for which  $x\leq y.$ 

In particular, one may consider designs on the hyperoctahedron. For instance, in $\Hy^{3}$ there is a set $B$ of three vertex-disjoint $2$-faces (edges) which partition the six vertices of $\Hy^{3}.$ It follows that $B$ is a $(1,2,1)$-design on $\Hy^{3}.$ Similarly, if $x$ is an $n$-set in $\Hy^{n}$ then $B=\{x,\,\bar{x}\}$ is a $(1,n,1)$-design on $\Hy^{n}$, and so on. 

\medskip 
\begin{thm}\label{design} Let $0<k\leq n$ and let  $B$ be a family of $k$-faces of $\Hy^{n}.$ Put $f:=\sum\,x$ (with $x\in B$) and let $f_{k,j,i}\in E_{k,j,i}$ denote the spectral components of $f.$  Then $B$ is a $(t,k,\ell)$-design on $\Hy^{n}$ for some $0<t\leq k$ and  $\ell$ (depending on $t)$ if and only if \\[10pt]
(a) $$f_{k,0,0}=\frac{\ell}{2^{k-t}{n-t\choose k-t}}\,\,l_{k}\quad\text{with}\quad l_{k}=\sum_{x\in L^{n}_{k}}\,x \,\,\text{,\,\,\,and}$$\\[10pt]
(b) \,\, for all $(j,i)$ with  $0\leq j\leq t,$ $0\leq i\leq \min\{t-j,\, n-k\}$ and $(j,i)\neq (0,0)$ we have $$f_{k,j,i}=0\,.$$
\end{thm}

\bigskip 
Since the components of $f=f_{B}$ can be computed explicitly using Proposition~\ref{compo}\, we obtain one equation for $f$ for every $(j,i)$ with  $0\leq j\leq t$ and  $0\leq i\leq \min\{t-j,\, n-k\}$ in terms of $\gnp_{k,j}$ and $\s_{k}.$ These include Fisher type inequalities for the number of blocks etc, for which we however omit details here. 

The precise analogue of Theorem~\ref{design}\, for ordinary $t$-designs (on the Boolean Lattice) is due to Graver and Jurkat~\cite{GJ}. It can be proved by the same arguments that follows now.

\bigskip
\pf  Suppose that  $B$ is a $(t,k,\ell)$-design and $f=\sum_{x\in B}\,x.$ By counting pairs $x\subseteq y$ with $x\in L^{n}_{t}$ and $y\in B$ we have  $|L^{n}_{t}|\cdot \ell=|B|\cdot {k\choose t}$ and hence by \,(\ref{deco02})\,
$$f_{k,0,0}=\frac{|L^{n}_{t}|\ell}{|L^{n}_{k}|{k\choose t}}\,l_{k}=\frac{\ell}{2^{k-t}{n-t\choose k-t}}\,\,l_{k}\,.$$

Further, if $x$ is a $t$-face then $\e_{t}^{k}(x)$ is the sum over all $k$-faces containing $x.$ Therefore $(\e_{t}^{k}(x),f)=\ell$ and so  $(\e_{t}^{k}(x-x'),f)=0$ for any two $t$-faces $x$ and $x'.$ The space $U$ generated by all $x-x'$ with $x,\,x'\in L^{n}_{t}$ has dimension $|L^{k}_{t}|-1$ and is perpendicular to $l_{t}.$ Hence $W=\oplus E_{t,j,i}$ where the sum runs over all $0\leq j\leq t,$ $0\leq i\leq \min\{t-j,\, n-t\}$ and $(j,i)\neq (0,0).$ It follows from part (iii) of Theorem~\ref{fulldecomp} that  $\e_{t}^{k}(W)=\oplus E_{k,j,i}$ where the sum runs over all $0\leq j\leq t,$ $0\leq i\leq \min\{t-j,\, n-t,\,k-j,\, n-k\}=\min\{t-j,\, n-k\}$ and $(j,i)\neq (0,0).$ Hence $f_{k,j,i}=0$ when the indices are in this range. 

Conversely, if $f_{k,j,i}=0$ for all $0\leq j\leq t,$ $0\leq i\leq \min\{t-j,\, n-k\}$ and $(j,i)\neq (0,0)$ then by the same argument we have $(\e_{t}^{k}(x-x'),f)=0$ for any two $t$-faces $x$ and $x'.$ Hence $\ell^{*}:=(\e_{t}^{k}(x),f)$ is independent of $x$  and  by \,(\ref{deco02})\, we have that $B$ is a $(t,k,\ell)$-design with $\ell^{*}=\ell.$ \dne

\medskip
{\sc 3.4 Association Schemes:\,} We mention the connection to the association schemes on $\Hy^{n}.$ Fix some $0\le k\leq n$ and let $0\leq j\leq k.$ Then two $k$-faces $x$ and $x'$ of $\Hy^{n}$ are {\it $j$-associates,} denoted   $x\,\sim_{j}x',$   if $|x\cap x'|=k-j.$ It is well-known that this relation defines an association scheme on the set of all $k$-faces on $\Hy^{n},$ see Chapter 11 of \cite{BI}  and \cite{JG}.

Let $\rho_{j}\!:\,M_{k}^{n}\to M_{k}^{n}$ be the corresponding association  map, given by  $$\rho_{j}(x)=\sum_{x\,\sim_{j}\,x'}\,x'.$$ 
For instance, $\rho_{0}(x)=x$ and $\rho_{1}(x)=\gnm_{k}(x)- kx.$ More generally, it is easy to show that $\rho_{j}$ is a polynomial in $\gnm_{k}$ of degree $j$ for any $j\leq k.$ Therefore all parameters, eigenspaces and idempotents of the association scheme are determined by Theorem~\ref{fulldecomp} and Corollary~\ref{cor2.8}. It is also clear that the incidence maps $\p$ and $\e$ are the appropriate tools to compare such schemes, as $k$ varies, at the level of their associated modules. While there is the notion of a homomorphism from one association scheme to another,  this is not the suitable tool for this particular purpose, see Zieschang~\cite{Z} and Xu~\cite{X}. (Indeed, the same is true for ordinary permutation modules: many essential features can not be described at the level of permutation sets only.) From this point of view the spectral decompositions constructed here can be obtained for {\it any} ranked partially ordered set from its incidence maps, independently of whether or not the set support an  association scheme. 

\bigskip\bigskip

\section{\sc Group Actions}

In this section we consider automorphisms of $\Hy^{n}$ and the associated permutation modules. As before let $V=\{\ga_{1}, \gab_{1},..,\ga_{n}, \gab_{n}\}$ be a set of $2n$ paired vertices with $\bar{\gab}_{i}=\ga_{i}$ for all $1\leq i\leq n.$ 
Then  the hyperoctahedral group $B_{n}$ is the group of all permutations $g$ in ${\rm Sym}(V)$  such that $\gbb^{g} = \overline{\gb^{g}}$ for all $\gb \in V.$ 
The action on $V$ can be extended to the faces of $\Hy^{n}$  by setting $$g\!:x=\{\gb_{1},..,\gb_{k}\}\,\,\to\,\, x^{g}=\{\gb_{1}^{g},..,\gb_{k}^{g}\}$$ for any face $x\in L^{n}_{k}.$ It is clear that $B_{n}$ preserves the partial order on $L^{n}$ and permutes $L^{n}_{k}$ transitively for each $k.$ (It  is not needed here but we mention that $B_{n}$ is the full automorphism group of  ${\mathcal H}^{n}.)$ In particular, if $\F$ is a field then $\F L^{n}_{k}=M^{n}_{k}$ is a transitive $\F B_{n}$ permutation module for each $k.$ Throughout we suppose that the characteristic  of $\F$ is $p>n$ or $p=0.$

\bigskip
{\sc 4.1 Irreducibility:\,} The elements of  $B_{n}$ commute with all the maps $M^{n}\to M^{n}$ mentioned in the previous sections. Therefore the eigenspaces $E_{k,j}$ and $E_{k,j,i}$ defined in Section 2 are invariant under $B_{n}.$ The ring of all linear maps $M^{n}\to M^{n}$ which commute with $B_{n}$ is denoted by ${\rm Hom}_{B_{n}}(M^{n}).$

\newpage
\begin{thm}\label{thm3.1} Let $\F$ be a field of characteristic $p>n$ or $p=0.$ Let $0\leq  k\leq n.$ Then the following hold: \\\\[-10pt]
(i) \,\, Let $0\leq j\leq k$ and $0\leq i \leq \min\{k-j, \,n-k\}.$  Then $E_{k,j,i}$ is an irreducible $\F B_{n}$-module. \\\\[-10pt]
(ii) $M_{k}$ is a multiplicity free $\F B_{n}$-module. The restrictions $\e_{k}$ and $\p_{k}$ with $0\leq k\leq n$ generate ${\rm Hom}_{B_{n}}(M^{n})$ as a ring. 
\end{thm}

For $m=\min\{k-j,n-k\}$ the modules $E_{k,j,m}$ along  the bottom diagonal in {\sc Figure}~2 are the Specht modules of $B_{n}$ for certain double $2$-part partitions in the standard representation theory of $B_{n},$ see   \cite{DJ, GeKi}.

\medskip
\pf {\it(i)} \,Let $x=\{\ga_{1},..., \ga_{k}\}\subseteq {V}.$ Since $G:=B_{n}$ is transitive on $L^{n}_{k}$ the permutation rank of $G$ on $L^{n}_{k}$ is the number $r$ of $H$-orbits on $L^{n}_{k}$ when $H$ is the setwise-wise stabilizer of $x.$ As a permutation group on $V$ 
we have $H={\rm Sym}(k)\times B_{n-k}$ where the second group is the octahedral group on $V\setminus (x\cup\bar{x})=V'$ and where the first group has two repeated actions \,\big(the diagonal action of ${\rm Sym}(k)\times {\rm Sym}(k)\big)$\,  of ${\rm Sym}(k)$ on $\{\ga_{1},..., \ga_{k}\}$ and $\{\gab_{1},..., \gab_{k}\}.$ It is clear that each $H$-orbit $y^{H}$ is characterized by the two numbers 
$a:=|y\cap x|$ and $b:=|y\cap V'|.$ Therefore  $r$ is the number of pairs $(a,b)$ with $0\leq a\leq k$ and $0\leq b\leq \min\{k-a,\,n-k\}.$ (The remaining $k-a-b$ points of $y$ can be chosen in $\bar{x}.)$ 

The permutation rank of a transitive group can be obtained also in a different way. Without loss assume that $\F$ is algebraically closed. Let $I_{1},\,I_{2},\,...,\,I_{t}$ be the irreducible characters of $B_{n}$ over $\F$ and let $n_{s}$ be the multiplicity of $I_{s}$ in $M^{n}_{k}.$ (Here we need that the characteristic  of $\F$ is $p>n$ or $p=0.)$ Then the permutation rank satisfies $r=\sum_{s=1}^{t}\,n_{s}^{2}.$ According to Theorem~\ref{fulldecomp} the number $R$ of summands  in $M_{k}=\bigoplus_{j,i} E_{k,j,i}$ is equal to the number of pairs $(j,i)$ with $0\leq j\leq k$ and $0\leq i\leq \min\{k-j,\,n-k\}.$  Comparing this to the count for $r$ above we have $R=r.$ Since $R$ is a lower bound for $r=\sum_{s=1}^{t}\,n_{s}$ we have $r=\sum_{s=1}^{t}\,n_{s}=\sum_{s=1}^{t}\,n_{s}^{2}.$ Hence $n_{s}\in \{0,1\}$ for all $s$ and so  $E_{k,j,i}$ is irreducible for all  $0\leq j\leq k$ and $0\leq i \leq \min\{k-j, \,n-k\}.$ 

{\it (ii)}\, It also follows that  $E_{k,j,i}$ has multiplicity $1$  in $M_{k}$ for  $0\leq j\leq k$ and $0\leq i \leq \min\{k-j, \,n-k\}.$ Furthermore, by Proposition~\ref{compo}\, the projection maps $\pi_{i,j,k}$ with $0\leq j\leq k\leq n$ and $0\leq i\leq \min\{k-j, \,n-k\}$ can be computed from the restrictions $\e_{k}$ and $\p_{k}.$
The remainder follows by standard arguments from representation theory and the third part  of Theorem~\ref{fulldecomp}.
\dne

\bigskip
{\sc 4.2 Orbits of Subgroups:} \, Let $G$ be a subgroup of $B_{n}$ and let $0\leq k\leq n.$ Then $G$ permutes the $k$-faces of ${\mathcal H}^{n}$ and so we let  $d^{G}_{k}$ denote the number of $G$-orbits on $L^{n}_{k}.$ Orbit numbers play an important role in many investigations. For instance, one may ask when the transitivity of $G$ on $L^{n}_{k}$ implies transitivity on $L^{n}_{k-1}.$    

\medskip
As an example, consider the octahedron in dimension $3.$ Its  full rotation group is the subgroup ${\rm Sym}_{4}\subseteq B_{3}.$ It is transitive on $L^{3}_{0},\,L^{3}_{1},\,L^{2}_{2}$ and $L^{3}_{3}.$ On the other hand, ${\rm Alt}_{4}$ is transitive on $L^{3}_{0},\,L^{3}_{1}$ and $L^{3}_{3}$ but has two orbits on $L^{n}_{2}.$ These correspond to the two tetrahedra that are inscribed in the cube. See again {\sc Comment}~2\, following Theorem~\ref{rank}. The same situation arises in any dimension. The situation is genuinely different from the  rank-symmetric unimodal orbit numbers that occur for A-type geometries.

\medskip
So, as before assume that $\F$ has characteristic $p>n$ or $p=0.$ If $N$ is a  $G$-invariant submodule of $\F M^{n}$ let $$N^{G}:=\{\,f\in N\,:\, f^{g}=f \,\,\text{for all}\,\,  g\in G\,\}$$ be the submodule of all elements in $N$ that are fixed by $G.$ In other words, if $f=\sum\,\,f_{x}\,x$\, is in $N$ then $f$ \,is in $N^{G}$ if and only if $x,\,x'$ being in the  same $G$-orbit on $L^{n}$ implies that $f_{x}=f_{x'}.$ In the case $N=M^{n}_{k}$ denote $M^{G}_{k}:=(M^{n}_{k})^{G}.$

Let $O_{1},\,..,\,O_{d}$ with $d=d_{k}^{G}$ be the distinct $G$-orbits on $L^{n}_{k}.$ For each $1\leq \ell\leq d$ we put  $f_{\ell}:=\sum_{x\in O_{\ell}}\,x.$  From above it is clear  that $f_{1},\,..,\, f_{d}$ forms a basis of $M^{G}_{k}.$ In particular, $$d_{k}^{G}=\dim(M^{G}_{k}).$$
We call $M^{G}_{k}$ the {\it $k$-orbit module} of $G.$ It is clear that the $G$-orbits on $L^{n}_{k}$ explicitly is equivalent to having complete information about  $M^{G}_{k}.$ The Decomposition Theorem therefore translates into a result on orbits. Using the same notation as in Theorem~\ref{fulldecomp} we have the following corollary:

\bigskip

\bigskip
\begin{thm} \label{Orb2}Assume that $\F$ is a field of characteristic $p>n$ or $p=0.$ Let $G\subseteq B_{n}$ and let $0\leq  k\leq n.$ Then the following hold:\\[10 pt]
 (i)\,\,  $M^{G}_{k}=\bigoplus _{i,j} \,\,E_{k,j,i}^{G}$ with $0\leq j\leq k$ and $0\leq i\leq \min\{k-j,\,n-k\}.$ \\[10 pt] (ii) \, For every  $0\leq i\leq n-j$ the maps $\e$ and $\p$  restrict to isomorphisms $$\e, \p\!:\,E_{j+i,j,i}^{G}\,\,\longleftrightarrow\,\, E_{j+i+1,j,i}^{G}\,\,\longleftrightarrow\,\, \dots \,\,\longleftrightarrow\,\,  E_{n-i,j,i}^{G}$$  while $\p(E_{j+i,j,i}^{G})\,=\,0\,=\,\e(E_{n-i,j,i}^{G}).$ \end{thm}

Note that the $E_{k,j,i}^{G}$ can be computed explicitly from $M^{G}_{k}$ by use of  Proposition~\ref{compo}.

\pf The two  parts follow from Theorem~\ref{fulldecomp}:\, For (i) it suffices to note that if $K$ and $N$ are $G$-invariant submodules of $M^{n}$ with $K\cap N=\{0\}$ then $(K\oplus N)^{G}=K^{G}\oplus N^{G}.$ For (ii) it suffices to note that if $A$ is a $G$-module and if $\varphi\!:\,A\to K$  is a $G$-homomorphism then  $\varphi$ restricts to a $G$-homomorphism $G\!:\,A^{G}\to K^{G}.$  \dne 

A second corollary gives inequalities for the orbit numbers. 

\bigskip
\begin{thm} Let $G\subseteq B_{n}$ and let $2k\leq n+1.$   Denote the number of $G$-orbits on $L^{n}_{k}$ by $d^{G}_{k}.$ Then $$d^{G}_{k-1}\leq d^{G}_{k}.$$ 
Furthermore, if $O_{1},\,..,\,O_{d}$ with $d=d_{k}^{G}$ are the distinct $G$-orbits on $L^{n}_{k}$ and if $f_{\ell}:=\sum_{x\in O_{\ell}}\,x$ for $1\leq \ell\leq d$ denote the corresponding orbit sums then $d^{G}_{k-1}= d^{G}_{k}$ if and only if $\pi_{k,k,0}(f_{\ell})=0$ for all $0\leq \ell\leq d.$ In the case of equality $M^{G}_{k}$ and $M^{G}_{k-1}$ are isomorphic to each other as $G$-modules. 
 \end{thm}

\pf Using Theorem~\ref{Orb2} twice, together with its second part, we have that $M^{G}_{k}$ is isomorphic to $M^{G}_{k-1}\oplus E_{k,k,0}^{G}.$ Now use the fact that $E_{k,k,0}=\pi_{k,k,0}(M^{n}_{k}).$ \dne

\end{document}